\newif\ifarxiv%
\arxivtrue

\ifarxiv%
    \documentclass[OJMO,e-only]{cedram}
    \def\copyrightline{}
\else%
    \documentclass[OJMO,manuscript]{cedram}
\fi%

\usepackage[all]{xypic}
\usepackage{booktabs}

\newcommand{\R}{\mathbb{R}}

\newcommand{\Np}{\mathbb{N}_+}

\newcommand{\sphere}[1]{\mathbb{S}^{#1}}

\newcommand{\con}{\mathbf{c}}
\DeclareMathOperator{\argmin}{arg\,min}

\newcommand{\vspan}{\mathrm{span}} % Span of vectors
\newcommand{\vect}{\mathrm{vec}} % Vectorisation
\newcommand{\unvect}{\mathrm{unvec}} % De-vectorisation

\newcommand{\E}{\mathbb{E}} % Expectation
\newcommand{\distr}{\sim} % Distributed according to
\newcommand{\iiddistr}{\overset{\text{i.i.d.}}{\distr}} % iid 

\newtheorem{assumption}[theo]{Assumption}

\usepackage[ruled,vlined]{algorithm2e}
\title[Linearly Separable CBO]{Consensus-based optimization for linearly separable functions}
\author{\firstname{Christian} \lastname{Fiedler}}
\address{TUM School of Computation, Information and Technology
and Munich Center for Machine Learning (MCML) \\
Boltzmannstr. 3, 85748 Garching bei M\"unchen,
Germany}
\email{christian.fiedler@tum.de}
\author{\firstname{Tim} \lastname{Roith}}
\address{TUM School of Computation, Information and Technology and Munich Center for Machine Learning (MCML)\\
Boltzmannstr. 3, 85748 Garching bei M\"unchen, Germany}
\email{tim.roith@tum.de}
\thanks{CF and TR acknowledge the support of the Munich Center for Machine Learning and the ERC Advanced Grant NEITALG, grant agreement No. 101198055. %
Funded by the European Union. 
{\nolinenumbers\\[.4em]}
\begin{minipage}{.3\textwidth}
\includegraphics[height=1.5cm]{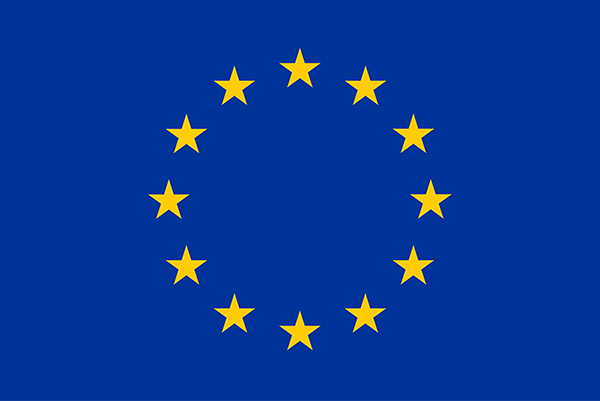}\hspace{1em}%
\includegraphics[height=1.5cm]{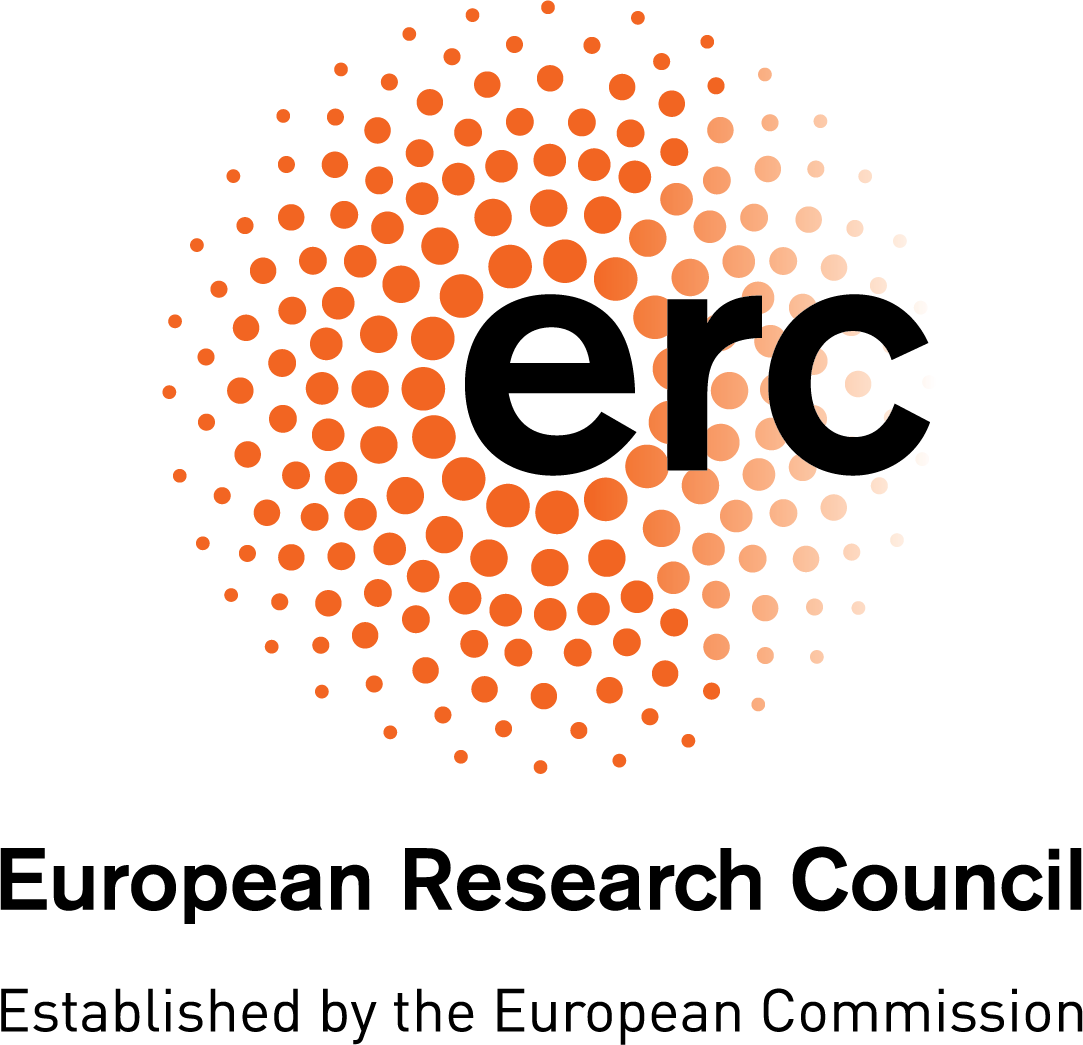}
\end{minipage}%
\begin{minipage}{.7\textwidth}
Views and opinions expressed are however those of the authors only and do not necessarily reflect those of the European Union or the European Research Council Executive Agency. 
Neither the European Union nor the granting authority can be held responsible for them. This project has received funding from the European Research Council (ERC) under the European Union's Horizon Europe research and innovation programme (grant agreement No. 101198055, project acronym NEITALG).
\end{minipage}%
{\nolinenumbers\\[.4em]}
In accordance with the Leiden declaration on AI and mathematics we disclose that artificial intelligence has been used in the preparation of this manuscript. Concretely Claude Opus 4.8 was employed to assist with the code for the numerical experiments and proofreading. The manuscript is written completely by the authors and we take full responsibility for the content.%
}%
\begin{abstract}%
Consensus-based optimization (CBO) is an efficient metaheuristic for global optimisation with attractive mathematical properties, allowing global convergence results even in non-convex settings. 
In practice it suffers greatly from the curse of dimensionality, as do most particle-based optimisers. 
Different strategies have been proposed to apply CBO even for high-dimensional optimisation problems, the most prominent being the so-called anisotropic noise model. 
However, a recent work \cite{bonandin2026exploiting} highlighted that this method performs well primarily on separable objective functions, and even small coordinate rotations severely degrade the performance.
Motivated by this observation, in this work we study the case where the objective function is separable only in a linearly transformed coordinate system. 
To leverage the efficiency of anisotropic CBO we propose a reparametrization scheme, which provably estimates such a coordinate transformation and then applies CBO in the new variables. 
Furthermore, we show how our scheme can be interpreted as noise adaptation in CBO.
Numerical examples highlight the efficacy of the method, demonstrating significant performance improvements on challenging benchmark objectives.

\end{abstract}
\begin{document}%
\maketitle%

\section{Introduction}
Consider the general global optimisation problem
\begin{align} \label{eq:globalOptProblem}
\min_{x\in X} f(x)
\end{align}
with objective function $f:X\to\R$ and search or decision space $X\subseteq \R^d$. We are particularly interested in the case of non-convex $f$ accessible only through noise-free queries, i.e., a classic zeroth-order optimisation problem with oracle function access.
In particular, we do not assume access to derivatives of $f$. 
A promising method for this setting is consensus-based optimization (CBO) \cite{pinnau2017consensus}, a meta-heuristic based on an interacting particle system.
CBO comes with strong theoretical guarantees \cite{fornasier2024consensus,fornasier2022convergence,carrillo2021consensus,bonandin2025strong,riedl2024mathematical} and achieves good empirical performance \cite{doi:10.1142/S0218202525500563,beddrich2026constrained}, with applications in machine learning \cite{carrillo2021consensus,roith2025consensus,fornasier2021consensus}, PDEs \cite{khatab2026consensus,beddrich2026constrained}, inverse problems \cite{bouillon2026localized,doi:10.1142/S0218202525500563} and robotics \cite{sun2026consensus}. To deal with high-dimensional problems of the form \eqref{eq:globalOptProblem}, a variant of CBO has been introduced that uses anisotropic noise \cite{carrillo2021consensus}, cf. Section \ref{sec:cbo} for details.
Anisotropic CBO also enjoys strong theoretical guarantees \cite{fornasier2022convergence} and exhibits good empirical performance even on challenging high-dimensional problems.
However, it was recently observed in \cite{bonandin2026exploiting} that anisotropic CBO performs well primarily for separable objective functions,
i.e., functions of the form
\begin{equation}
\label{eq:separableFunction}
f(x) = \sum_{j=1}^d f_j(x_j), \qquad f_j:\R\to\R,\text{ for } j=1,\ldots,d.
\end{equation}
Geometrically, this means that the objective function varies only along the coordinate axes,
and it is a property extensively studied and used in optimisation, statistics and machine learning (where such functions are often called additive).
It was observed in \cite{bonandin2026exploiting} that even slight deviations from this property severely degrade performance of anisotropic CBO.
For example, let $f$ be of the form \eqref{eq:separableFunction} and $Q$ describing a small-angle rotation.
For many cases of $f$, including prominent benchmark functions, anisotropic CBO performs very well on $f$,
but very badly on $\tilde{f}(x)=f(Qx)$.
We reproduce and illustrate this effect in Figure~\ref{fig:rotrast}. 
Motivated by this problem, in this work we consider functions of the form
% In most cases we are interested in the global optimization case where $X=\R^d$. 
% We are especially interested in situations, where $f$ is not convex and we cannot access its derivatives. However, the specific assumption that we make within this work, is that $f$ is additive in the following sense,
%
\begin{equation} \label{eq:linearlySeparableFunctions}
f(x) = \sum_{j=1}^m f_j(\langle x, w_j \rangle) %=: g(Wx), \quad \text{ where } g:X\to\R^d \text{ is separable and } W=(w_1,\ldots, w_m)
\end{equation}
with $1 \leq m \leq d$ and $w_1,\ldots,w_m \in \R^d$ are linearly independent.
%The example above is a special case with $m=d$ and $w_1,\ldots,w_m$ given by the rows of the rotation matrix $Q$.
We can interpret \eqref{eq:linearlySeparableFunctions} as a neural network with one hidden layer and a linear activation function in the output layer, providing a bridge to \cite{fornasier2021robust,fornasier2025efficient} on which we build later on.
Furthermore, defining $g: \R^m \rightarrow \R$, $g(z)=\sum_{j=1}^m f_j(z_j)$, and $W=\begin{pmatrix} w_1 & \cdots & w_m \end{pmatrix}^\top$,
we can write $f(x)=g(Wx)$ and hence interpret \eqref{eq:linearlySeparableFunctions} also as a special case of the multi-index model.
The example above corresponds to the case $m=d$ and $W$ given by the rotation matrix $Q$.
% Note that \eqref{eq:linearlySeparableFunctions} can be interpreted as a neural network with one hidden layer and a linear activation function in the output layer. In our setting we have noise-free access to function values at arbitrary inputs.
% This crucial observation allows us to use the tools from \cite{fornasier2021robust,fornasier2022robust,fiedler2019learning}.

The importance of functions of the form \eqref{eq:linearlySeparableFunctions} is two-fold.
First, intuitively functions of this form are separable in an appropriate coordinate system.
% More precisely, given a function $f$ as in \eqref{eq:linearlySeparableFunctions}, if we define $W=\begin{pmatrix} w_1 & \cdots w_m\end{pmatrix}$ and $\tilde f: \R^m \rightarrow \R$, $\tilde f(z)=f(Wz)$ (ignoring domain issue for simplicity), then $\tilde f$ is separable, i.e., of the form \eqref{eq:separableFunction}.
For this reason, we call functions \eqref{eq:linearlySeparableFunctions} also \emph{linearly separable}.
From a practical point of view this fact is important since we might not know the coordinate system in which a function is separable.
Assuming just \eqref{eq:linearlySeparableFunctions} instead of \eqref{eq:separableFunction} is therefore more realistic.
Second, functions of the form \eqref{eq:linearlySeparableFunctions} can be effectively much lower-dimensional than the ambient dimension $d$, namely if $m \ll d$, in which case also optimisation problems such as \eqref{eq:globalOptProblem} are effectively lower dimensional. 
\paragraph*{Contribution}
In this work, we propose a global optimisation approach that takes advantage of the special structure in \eqref{eq:linearlySeparableFunctions}.
%Optionally, we start with a dimensionality reduction step (in case $m \ll d$), for which we use the method from \cite{fornasier2021robust}.
We use the efficient method from \cite{fornasier2021robust,fornasier2022robust,fiedler2023stable} to estimate the vectors $w_1,\ldots,w_m$ purely from function evaluations.
In turn, these estimates are used in a CBO scheme to reduce the global optimisation problem \eqref{eq:globalOptProblem} to one involving an approximately separable objective function. We consider three options, namely reparametrisation, noise modification, and componentwise splitting. We demonstrate the effectiveness of the method on different benchmark functions and compare its performance to CMA-ES.
%
% where $w_1,\ldots, w_m\in \R^{d}$ are unknown vectors and $\phi:\R\to\R$ is some function. Furthermore, we denote by $W=(w_1|\ldots|w_m)\in\R^{m\times d}$ the matrix composed of these vectors, and in the special case where $d=m$ and $W$ is an orthogonal matrix, this then means that
% %
% \begin{align*}
% \tilde{f}(x) := f(Wx) = \sum_{i=1}^m \phi\langle w_i, Wx\rangle) = 
% \sum_{i=1}^m \phi(x_i)
% \end{align*}
% %
% is a separable function.
%
%
% \paragraph{Outline} In Section \ref{sec:background}, we collect background and preliminaries.
% In particular, we introduce (anisotropic) CBO and its convergence guarantees in Section \ref{sec:cbo},
% the dimensionality reduction approach in Section \ref{sec:dr},
% and the vector estimation method in Section \ref{sec:hessianSampling}.
% We describe our proposed approach in Section \ref{sec:method} and providing theoretical guarantees,
% and in Section \ref{sec:experiments} we demonstrate the efficacy of our method in numerical experiments.
% We close in Section \ref{sec:conclusion} with a summary and discussion.
 
\section{Background and preliminaries} \label{sec:background}
\subsection{Consensus-based optimisation} \label{sec:cbo}

Consensus-based optimisation was introduced in \cite{pinnau2017consensus} as a particle-based optimisation scheme with a well-posed mean-field limit. It employs an ensemble of particles $\mathbf{x} = (x^{(1)}, \ldots, x^{(n)})$ exploring the state space $X$ through diffusion and drifting towards the consensus point
\begin{align}\label{eq:con}
\con_\alpha^f(\mathbf{x}) := 
\sum_{i=1}^n \omega_i^f \cdot x^{(i)},\qquad 
\omega_i^f:= \frac{\exp(-\alpha f(x^{(i)}))}{\sum_{j=1}^n \exp(-\alpha f(x^{(j)}))},
\end{align}
which is a weighted average, where $\alpha>0$ controls the weighting strength. 
Given an initial distribution $\rho_0$ and 
%$x_{0}^{(i)}\sim \rho_0, i=1,\ldots,n$, i.i.d., 
$x_0^{(1)},\ldots,x_0^{(n)} \iiddistr \rho_0$,
the Euler--Maruyama discretised update step of the $i$-th particle then reads
\begin{align*}
x^{(i)}_{k+1} = 
x^{(i)}_k - \tau\, (x^{(i)}_k - \con^f_\alpha(\mathbf{x}_k)) + \sigma\, \sqrt{\tau}\, \mathsf{D}(x^{(i)}_k - \con^f_\alpha(\mathbf{x}_k))\, \xi_k^{(i)},
\end{align*}
where $\xi_k^{(i)}\sim \mathcal{N}(0, I_{d})$, $\tau>0$ denotes the step size and the matrix-valued function $\mathsf{D}:\R^d\to\R^{d\times d}$ determines the noise model. In particular, the following choices lead to isotropic and anisotropic noise models,
\begin{align*}
\mathsf{D}^{\text{iso}}(z) := \|z\|_2\cdot I_{d}, \qquad
\mathsf{D}^{\text{aniso}}(z) := \operatorname{diag}(z_1, \ldots, z_d),
\end{align*}
since in the first case the noise is scaled equally in every spatial direction, whereas in the second case each component receives a separate scale. At this point we want to highlight that this separation is tied to the standard basis, which is precisely the reason why anisotropic noise works especially well for separable functions and the extension to more general settings is the main motivation of this work.

\paragraph*{Convergence guarantees of CBO} 
CBO and its variants are supported by a rich convergence analysis which is an active area of research \cite{fornasier2024consensus,fornasier2022convergence,carrillo2021consensus,bonandin2025strong,riedl2024mathematical}. The main motivation for our study here is the recent work \cite{bonandin2026exploiting}, which analyses how anisotropic diffusion in CBO exploits the structure of the underlying objective function. Assume that $g:\R^d\to\R$ is a fully separable function with unique global minimum $x^*$ and further satisfies suitable regularity and tractability assumptions (see \cite{bonandin2026exploiting} for the precise assumptions). Given an arbitrary accuracy level $\varepsilon>0$, \cite[Thm.~3.8, Cor.~3.9]{bonandin2026exploiting} show that with probability at least $(1-\delta_1)^d - \delta_2$, where $\delta_1\in (0,1-2^{-1/d}),\delta_2\in(0,1/2)$, one obtains for some $\vartheta\in (0,1)$
\begin{align*}
\left\|
\frac{1}{n} \sum_{i=1}^n x_{\lceil T/\tau\rceil}^{(i)} - x^*
\right\|_{\infty} \leq \varepsilon,\qquad \text{where}\quad T\geq \frac{C - \log(\varepsilon^2\delta_1)}{(1-\vartheta)(2-\sigma^2)} +1,\quad  n\geq \frac{8 C_{\text{MFA}}}{\varepsilon^2\delta_2},\quad 
\tau \leq \frac{\varepsilon^2\delta_2}{8C_{\text{NA}}}.
\end{align*}
Here, the constant $C$ depends on the initial distribution $\rho_0$, the constant $C_{\text{MFA}}$ is independent of $d$ and $C_{\text{NA}}$ depends only linearly on $d$. Running the CBO algorithm for a number of $K=\lceil T/\tau\rceil$ steps requires $K\cdot n$ number of queries of $g$. As noted in \cite[Rem.~3.10]{bonandin2026exploiting} the computational complexity to achieve the above error tolerance thus also scales linearly in the dimension $d$. Without assuming separability the constant $C_{\text{MFA}}$ would scale at least exponentially in $d$ and thus the computational complexity would do so, too. This highlights the immense computational advantage that anisotropic CBO gains from exploiting this special structure.

% \subsection{Dimensionality reduction through finite-differences} \label{sec:dr}
% Consider a linearly separable function of the form \eqref{eq:linearlySeparableFunctions} with $m\ll d$.
% In this case, it is beneficial to first estimate the $m$-dimensional subspace $\mathcal{V}=\vspan\{w_1,\ldots,w_m\}\subseteq\R^d$,
% since this allows to reduce the effectively dimensionality of the optimisation problem \eqref{eq:globalOptProblem} and the estimation procedure described in the next subsection.

\subsection{Efficient separability recovery through Hessian sampling}\label{sec:hessianSampling}
For ease of exposition, we assume that $X=\R^d$ and that $f_1,\ldots,f_m$ are twice continuously differentiable.
Without loss of generality we assume that $\|w_1\|_2=\ldots=\|w_m\|_2=1$, since $w_i\not=0$ by linear independence and hence we can rewrite $f_j(\langle x, w_j\rangle)=f_j(\|w_j\|_2 \cdot \langle x, \|w_j\|_2^{-1} w_j \rangle)$. Denote by $H[f](x)$ the Hessian of $f$ at $x$ and observe that
\begin{equation}
    H[f](x) = \sum_{i=1}^m f_i''(\langle x, w_i\rangle) w_i \otimes w_i 
    \: \in \: 
    \vspan\{ w_1 \otimes w_1,\ldots, w_m \otimes w_m \} =: \mathcal{W},
\end{equation}
where we defined $v \otimes w = vw^\top$ for vectors $v,w\in\R^d$.
It was discovered in \cite{fornasier2021robust} that $\mathcal{W}$ can be approximated with high probability with random queries and numerical differentiation of $f$.
In turn, the individual $w_i$'s can be recovered from the approximated matrix subspace with a randomised algorithm.
We will describe this method now in detail: For $\epsilon>0$ let $\hat{H}_\epsilon[\cdot]$ be a numerical approximation method for the Hessian, e.g., finite differences, fulfilling the following assumption for some set $\emptyset\neq S\subset X$.
\begin{assumption} \label{assump:hessianApprox}
There exist constants $C_{\hat{H}},\epsilon_0\geq 0$ such that for all $\epsilon>\epsilon_0$ we have
\begin{equation}
    \| H[f](x) - \hat{H}_\epsilon[f](x)\|_F \leq C_{\hat{H}} m \epsilon \qquad \forall x \in S.
\end{equation}
\end{assumption}
For example, if $S$ is bounded, the $f_i$'s are three times continuously differentiable, and we use finite differences for $\hat{H}_\epsilon$, then Assumption \ref{assump:hessianApprox} holds with $\epsilon_0=0$,
cf. \cite[Lemma~4.1]{fornasier2021robust}.
We allow the case $\epsilon_0>0$ since computing $\hat{H}_\epsilon[f]$ might involve additional approximations that lead to an unavoidable error, e.g., if an initial dimension reduction has been applied. For a matrix $M$, denote by $\vect[M]$ the vector that results from stacking the columns of $M$,
and denote by $\unvect[\cdot]$ the reverse (splitting a vector and using the subvectors as columns of a matrix). Let now $\mu$ be a Borel probability measure on $\R^d$ with support in $S$.
In practice, $S$ will usually be a centered ball and $\mu$ the uniform distribution on this ball.
\begin{assumption} \label{assump:identifiabilityHessianLinSepFunc}
The matrix
\begin{equation} \label{eq:identifiabilityMatrixWelldefined}
    \mathbb{H}_{\mu}[f] = \E_{X \distr \mu}\left[ \vect\big[ H[f](X) \big] \otimes \vect\big[ H[f](X) \big] \right]
\end{equation}
is well-defined and there exists a constant $\alpha_H >0$ such that the following holds:
\begin{equation} \label{eq:identifiabilityConditionHessianLinSep}
    \sigma_m(\mathbb{H}_{\mu}[f]) \geq \alpha_H\qquad\text{where } \sigma_m \text{ denotes the } m\text{th largest singular value}.
\end{equation}
\end{assumption}
The well-definedness of \eqref{eq:identifiabilityMatrixWelldefined} is a rather weak assumption,
and holds for example if $S$ is bounded.
The requirement \eqref{eq:identifiabilityConditionHessianLinSep} is essentially an identifiability condition and appeared for the first time in \cite[Theorem~4.2]{fornasier2021robust}.
\begin{assumption} \label{assump:hessianMaxSvalLinSepFunc}
There exists a constant $\bar{C}\geq 0$ such that for $X\distr \mu$ we have
\begin{equation}
    \|H[f](X)\|_F \leq \bar{C} m \qquad \text{almost surely}.
\end{equation}
\end{assumption}
The preceding condition is fulfilled under mild assumptions such as boundedness of $S$, cf. the proof of \cite[Theorem~4.2]{fornasier2021robust}.
If the vectors $w_1,\ldots,w_m$ fulfil an incoherence condition, one can derive a refined bound in this case as shown in the proof of \cite[Lemma~A.2]{fornasier2025efficient}. For $N_H\in\Np$, let $x^{(1)},\ldots,x^{(N_H)} \iiddistr \mu$, and define for $\epsilon>\epsilon_0$ the matrix
\begin{equation}\label{eq:hatY}
\hat{Y} = \begin{pmatrix}
\vect\left[\hat{H}_\epsilon[f](x^{(1)}) \right]
& \cdots & 
 \vect\left[\hat{H}_\epsilon[f](x^{(N_H)}) \right]
\end{pmatrix}
\in \R^{d^2 \times N_H}.
\end{equation}
Consider its SVD
$
\hat{Y} = \begin{pmatrix}
    \hat{U} & \hat{U}_0
\end{pmatrix}
\begin{pmatrix}
    \hat{\Sigma} & 0 \\
    0 & \hat{\Sigma}_0
\end{pmatrix}
\begin{pmatrix}
    \hat{V} & \hat{V}_0
\end{pmatrix}^\top
$
with $\hat{\Sigma}\in\R^{m \times m}$ and $\hat{\Sigma}_0 \in \R^{(d^2-m)\times(d^2-m)}$,
and the induced subspace
$
    \hat{\mathcal{W}} = \vspan\{ \unvect\left[ \hat{U}_{\cdot i} \right] \mid i = 1,\ldots, m \} \subseteq \R^{d \times d}.$ 
The next result, which follows by minor modifications of the proof of \cite[Theorem~4.2]{fornasier2021robust}, shows that $\hat{\mathcal{W}}$ is a good approximation of $\mathcal{W}$.
\begin{prop}\label{prop:a}
Under Assumptions \ref{assump:hessianApprox}, \ref{assump:identifiabilityHessianLinSepFunc}, and \ref{assump:hessianMaxSvalLinSepFunc},
for all $\epsilon > \epsilon_0$ and $s\in(0,1)$ with
\begin{equation}
    \sqrt{(1-s) \alpha_H} > C_{\hat{H}}  m\epsilon,
\end{equation}
denoting by $P_\mathcal{V}$ the orthogonal projector onto a subspace $\mathcal{V}$ and the operator norm by $\|\cdot\|_{F\to F}$ we have
\begin{equation}
     \| P_{\mathcal{W}} - P_{\hat{\mathcal{W}}} \|_{F \rightarrow F} 
     \leq 
    \frac{2 C_{\hat{H}} m\epsilon}{\sqrt{(1-s) \alpha_H} - C_{\hat{H}} m\epsilon}
\end{equation}
with probability at least
\begin{equation}
    1 - m\exp\left(-\frac{s^2 N_H \alpha_H}{2 \bar{C}^2 m^2} \right).
\end{equation}
\end{prop}
Note that for $\epsilon_0=0$, this result is applicable for any $s\in(0,1)$ by choosing $\epsilon>0$ small enough.
Furthermore, when using finite differences, then the sample complexity (number of function queries) will scale quadratically with $d$, but not exponentially.
Finally, in practice $m$ is determined by a gap in the singular value spectrum of $\hat{Y}$.
\begin{rema}
If $m<d$, then our approach in Section \ref{sec:method} also performs an implicit dimensionality reduction.
If $m \ll d$, then an initial dimensionality reduction can be useful due to the quadratic sample complexity of the Hessian approximation,
cf. \cite{fornasier2021robust} and \cite{fornasier2012learning} for such approaches.
\end{rema}
Consider now the nonconvex optimisation problem
\begin{equation}\label{eq:spmOpt}
    \max_{u \in \sphere{d-1}} \| P_{\hat{\mathcal{W}}}(u \otimes u) \|_F^2.
\end{equation}
One can derive it from geometric considerations, cf. \cite{fiedler2019learning}, and it is a special case of the subspace power method for tensor decomposition \cite{kileel2021landscape}.
Define now $\rho_2 = \sup_{x \in \sphere{d-1}} \sum_{j=1}^m \langle x, w_j \rangle^2 - 1,$
which is a measure of separation or incoherence of the vectors $w_1,\ldots,w_m$.
For example, for $m=d$ and orthonormal vectors, it is zero.
We now recall a result on the maximisers of \eqref{eq:spmOpt}. It corresponds to \cite[Theorem~7]{kileel2021landscape} for $n=2$ (matrices instead of general tensors) and $s=2$.
\begin{prop}\label{prop:b}
Define $\Delta_0 = \tfrac{1/6 - 10\rho_2}{22/3 - 20\rho_2}$
and
\begin{equation}
\hat{\mathcal{S}} = \left\{ x \in \sphere{d-1} 
\mid 
\| P_{\hat{\mathcal{W}}}(x \otimes x) \|_F^2 
\geq 
\frac{43-60\rho_2}{1-60\rho_2} \|P_{\mathcal{W}} - P_{\hat{\mathcal{W}}}\|_{F\rightarrow F}
\right\}.
\end{equation}
If $\|P_{\mathcal{W}} - P_{\hat{\mathcal{W}}}\|_{F\rightarrow F} < \Delta_0$,
then $\hat{\mathcal{S}}$ contains exactly $2m$ second-order critical points, each of which is a local maximiser of \eqref{eq:spmOpt},
and for each such point $\hat x$ there exists a unique $j \in \{1,\ldots,m\}$ and $\sigma\in\{-1,1\}$
with
\begin{equation}
\|\hat x - \sigma w_j \|_2^2 \leq \|P_{\mathcal{W}} - P_{\hat{\mathcal{W}}}\|_{F\rightarrow F}.
\end{equation}
\end{prop}
By finding all local maximisers described in the preceding result, we can therefore recover precise approximations of $w_1,\ldots,w_m$ with high probability.
To do so, we solve \eqref{eq:spmOpt} by sampling an initial value uniformly on the unit sphere and then run gradient ascent on the Riemannian manifold $\sphere{d-1}$ with step size $\gamma\geq 1/2$, cf. \cite{absil2008optimization} for an introduction to these optimisation methods.
The iterates converge linearly to a stationary point of the constrained optimisation problem, which is almost surely (w.r.t. the initialisation) a constrained local maximiser, cf. \cite[Theorem~9]{fiedler2023stable}.
To avoid spurious maximisers, one can simply check that the last iterate achieves a value close to 1, cf. \cite{fiedler2023stable,fornasier2025efficient} for practical considerations. We summarise the estimation algorithm in Algorithm \ref{alg:estW}.
\begin{algorithm}[H]\caption{Estimating $w_1,\ldots,w_m$ from queries to $f$} %$f = g\circ W$}
\label{alg:estW}
\DontPrintSemicolon%
\SetAlgoLined%
\SetAlgoNoEnd%
Sample $x^{(1)}, \ldots, x^{(N_H)}\iiddistr\mu$\;
Compute finite difference Hessians $\hat{H}_\epsilon[f](x^{(i)})$ for each $i=1,\ldots,N_H$, form $\hat{Y}$ in \eqref{eq:hatY}\;
Compute $\hat{\mathcal{W}}$ from SVD of $\hat{Y}$\;
\For{$k=1,\ldots, N_{\mathrm{restarts}}$}{%
\For{$l=1,\ldots, L$}{
$v^k_{l+1}:= \Pi(v_l^k + \gamma\nabla \|P_{\hat{\mathcal{W}}}(v_l^k\otimes v_l^k)\|_F^2)$, \quad where $v^k_0\sim \text{Uniform}(\sphere{d-1})$ and $\Pi(x):= x/\|x\|_2$%$V_l$
}
}%
Select $\{\tilde{w}_1,\ldots, \tilde{w}_{\tilde{m}}\}$ as the biggest linearly independent subset from $\{v_L^1, \ldots, v_L^{N_{\mathrm{restarts}}}\}$\;
\Return{$\tilde{w}_1,\ldots, \tilde{w}_{\tilde{m}}$}
\end{algorithm}

\section{Method} \label{sec:method}

In this section, we explain how we modify the CBO scheme to obtain an efficient algorithm for functions $f$ as in \eqref{eq:linearlySeparableFunctions}, i.e., which are separable in a transformed coordinate system. 
We collect the estimates $\tilde{w}_i$'s as rows of a matrix $\tilde{W}$, which is used as an approximation of $W$ (up to permutation and sign, which do not matter for the optimisation later on).
%We obtain an estimate $\tilde{W}$ of $W$
%as described in the previous section and then consider three different approaches that exploit the underlying structure of a linearly separable function. 
When the estimate is exact%, i.e., $\tilde{W}=W$, then 
we can apply the convergence guarantees from \cite{bonandin2026exploiting}. 
In our case Propositions~\ref{prop:a} and~\ref{prop:b} ensure a very good approximation with high probability,
and we leave the analytical treatment of the approximate case to future work.
%We leave the analytical treatment of the case where $\tilde{W}$ is only an approximation for future work.

\subsection{Reparametrising the objective function}

The most direct approach to employ the estimate $\tilde{W}$ is to reparametrise the objective function $f$. Denoting by $\tilde{W}^\dagger$ the pseudo-inverse of $\tilde{W}$, we define the objective
%
%\begin{align*}
$\tilde{g}(x) = f(\tilde{W}^\dagger x)$
%\end{align*}
%
which we can then assume to be (approximately) separable. 
Our benchmarks show that this methodology can greatly improve the performance of anisotropic CBO applied to $f$ itself.

\subsection{Noise modification through change of variables}

The previous reparametrisation can equivalently be shifted onto the noise term in CBO. 
Denote by $g$ the separable objective function such that $f(x)=g(Wx)$, where we assume $W\in\R^{d\times d}$ to be invertible. In its time-continuous form anisotropic CBO applied to $g$ yields the system of SDEs for $X_t=(X_t^{(1)}, \ldots, X_t^{(n)})$,
\begin{align*}
\mathrm{d} X_t^{(i)} = 
- (X_t^{(i)} - \mathbf{c}^g_\alpha(X_t))\mathrm{d}t + 
\sigma\, \mathsf{D}^{\text{aniso}}(X_t^{(i)} - \mathbf{c}^g_\alpha(X_t)) \, \mathrm{d}B_t^{(i)},\qquad \text{for } i=1,\ldots, n,
\end{align*}
where $B_t^{(i)}$ denote independent Brownian motions. With the change of variables $Y_t^{(i)} = W^{-1} X_t^{(i)}$, we see
\begin{align*}
\mathrm{d} Y_t^{(i)} = 
- (Y_t^{(i)} - \mathbf{c}^f_\alpha(Y_t))\mathrm{d}t + 
\sigma\,  W^{-1}\mathsf{D}^{\text{aniso}}(W(Y_t^{(i)} - \mathbf{c}_\alpha^f(Y_t))) \, \mathrm{d}B_t^{(i)},\qquad \text{for } i=1,\ldots, n,
\end{align*}
where we used the fact that 
\begin{align*}
W^{-1}\mathbf{c}^g_\alpha(X_t) 
=
\frac{\sum_{i=1}^n \exp(-\alpha g(X_t^{(i)})) W^{-1} X_t^{(i)}}{\sum_{i=1}^n \exp(-\alpha g(X_t^{(i)}))} = 
\frac{\sum_{i=1}^n \exp(-\alpha g(W Y_t^{(i)})) Y_t^{(i)}}{\sum_{i=1}^n \exp(-\alpha g(W Y_t^{(i)}))}=
\mathbf{c}^f_\alpha(Y_t).
\end{align*}
In other words, when we apply CBO to the transformed function $f$, we can employ the transformed noise model
\begin{align}\label{eq:noisetrafo}
\mathsf{D}^{\text{transform}}(z; W) := 
W^{-1} \mathsf{D}^{\text{aniso}}(Wz)
\end{align}
to restore the favourable convergence properties that CBO has for separable functions. 
%In the practical case, 
In numerical implementations, 
we then employ the noise model $\mathsf{D}^{\text{transform}}(\cdot; \tilde{W})$. 
This perspective connects to the covariance matrix adaptation evolution strategy (CMA-ES) \cite{hansen2023cmaevolutionstrategytutorial}, which is a related particle-based optimisation scheme,
% A possible connection to CBO is also pointed out in \cite{fornasier2026consensus}. 
cf. also \cite{fornasier2026consensus} for another possible connection to CBO.
CMA-ES adapts the covariance matrix that underlies the exploration online. Our experiments demonstrate how powerful a well-chosen noise model can be for CBO schemes, however, unlike CMA-ES it is a two-stage procedure. A combined adaptation scheme for CBO is left for future work.
\subsection{Splitting into one-dimensional problems}

%A further powerful technique for separable functions is componentwise splitting. 
Finally, we can also apply componentwise splitting.
%Namely, once 
Once we have obtained an approximately separable function $\tilde{g}$, we can consider $m$ one-dimensional optimisation problems
\begin{align*}
s_j := \argmin_{s\in\R} f(s \cdot \tilde{W}^\dagger e_j),\qquad\text{for } j=1,\ldots, m, \qquad \Rightarrow\quad  
x^* = \tilde{W}^\dagger (s_1,\ldots, s_m)^\top,
\end{align*}
where $e_j\in\R^{m}$ denotes the $j$-th standard basis vector. In our numerical examples we observe that using CBO for each sub-problem exhibits several practical advantages over running CBO on the full $d$-dimensional problem. 
Most prominently it is embarrassingly parallelizable with respect to the dimension $m$. Although now each single CBO run evaluates the function $f$ without sharing this information with other runs, the total number of required queries is typically lower than for the full problem, when both versions must reach the same accuracy level. 
\begin{rema}
For ease of notation let $m=d$, then \cite[Lem.~4.1]{bonandin2026exploiting} shows that for a measure $\rho=\rho^{(1)}\otimes\ldots\otimes \rho^{(d)}$
\begin{align*}
\con^g_\alpha(\rho) = \sum_{j=1}^d \con^{f_j}_\alpha(\rho^{(j)})\, e_j,\qquad\text{where}\quad 
\con_\alpha^h(\mu) := 
\int \omega_\mu^h(x)\cdot x \, \mathrm{d}\mu(x),\quad 
\omega_\mu^h(x):= \frac{\exp(-\alpha h(x))}{\int\exp(-\alpha h(x))\mathrm{d}\mu(x)},
\end{align*}
which would suggest that running the $d$ one-dimensional dynamics should be equivalent to the full $d$-dimensional dynamics. 
However, the identity above is only exact in the mean-field limit, so for finitely many particles the assumption $\rho=\rho^{(1)}\otimes\ldots\otimes \rho^{(d)}$ is typically not fulfilled.
The splitting approach thus indeed yields a different, and as we show in Section \ref{sec:experiments} often more efficient, numerical scheme.
\end{rema}

\section{Experiments} \label{sec:experiments}
\begin{figure}
\begin{subfigure}[t]{.23\textwidth}%
\includegraphics[width=\textwidth]{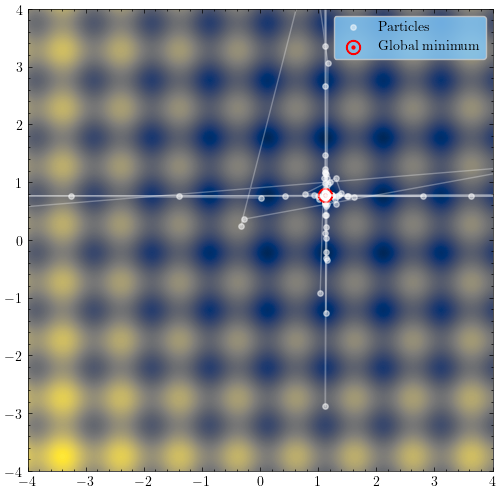}
\caption{Anisotropic CBO on the standard Rastrigin function.}
\end{subfigure}\hfill%
\begin{subfigure}[t]{.23\textwidth}%
\includegraphics[width=\textwidth]{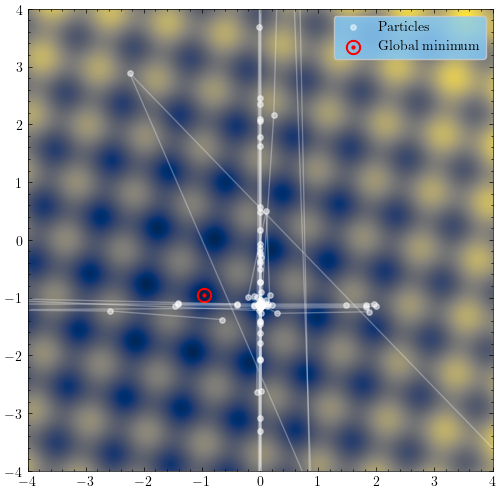}
\caption{Anisotropic CBO on the rotated Rastrigin function.}
\end{subfigure}\hfill%
\begin{subfigure}[t]{.23\textwidth}%
\includegraphics[width=\textwidth]{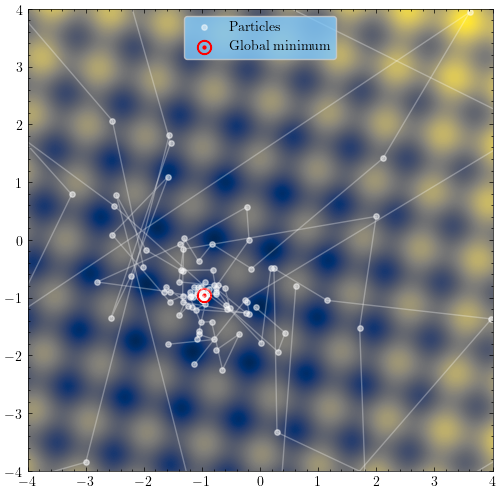}
\caption{Isotropic CBO on the rotated Rastrigin function.}
\end{subfigure}\hfill%
\begin{subfigure}[t]{.23\textwidth}%
\includegraphics[width=\textwidth]{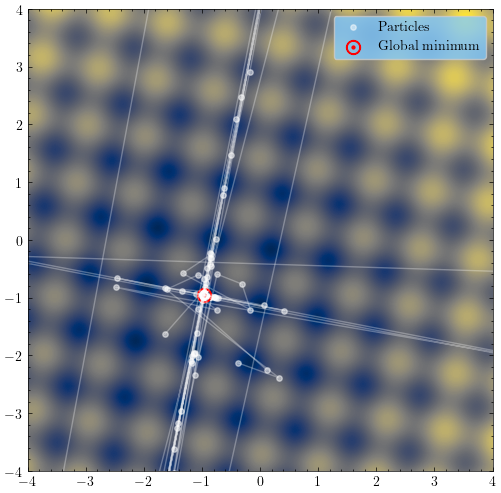}
\caption{CBO with the proposed rotated noise in \eqref{eq:noisetrafo}.}
\end{subfigure}\hfill%
\caption{We plot the evolution of the particles for different CBO algorithms. As discussed in \cite{bonandin2026exploiting} anisotropic CBO is highly efficient for separable functions. However, when we consider a rotated version, the algorithm still tries to exploit an underlying axis aligned structure, which does not match the structure of the objective. Isotropic CBO neither assumes nor exploits an underlying separability.}\label{fig:rotrast}
\end{figure}

We now study the numerical behaviour of our proposed method on different benchmark functions.
The code for the experiments is based on the \texttt{CBXPy} package \cite{Bailo_CBX_Python_and_2024} and available at \url{https://github.com/TimRoith/LinearlySeparableCBO}. 
%We test our pipeline on different benchmark functions, which we detail below. 
For our evaluation we allow each algorithm to employ a total number of $B=5{,}000{,}000$ queries to the objective function $f$. 
For the first stage of our method, we use the uniform distribution on the closed Euclidean ball of radius 20 as the sampling distribution $\mu$ and use $B_{H} = 2\, (d^2 + d)\,  N_H$ queries for the finite difference estimation of the Hessian, where in our experiments we choose $N_H=150$, denoting the number of Hessians to be sampled and thus $B_{H}$ uses around $60\%$ of the total budget. 
Note that the CBO schemes (which all employ $n=50$ particles) are then only allowed to use $B - B_{H}$ queries.
For the matrix $W$, we sample a random orthogonal matrix $Q\in\R^{d\times d}$ and a random invertible matrix $R\in\R^{d\times d}$, set $\hat{W}=Q+5 \cdot R$ and obtain $W$ via normalising the rows of $\hat{W}$.
A more detailed numerical study regarding optimal query allocation (for Hessian sampling and CBO) and the influence of possible matrix estimation errors is left for future work.
Finally, we also provide a comparison to the bi-population CMA-ES variant, see \cite{hansen2009benchmarking} using the implementation from the \texttt{pycma} package \cite{pycma}. 
We now detail the employed benchmark functions: The Rastrigin \cite{rastrigin1974systems} function is given by
\begin{align*}
g^{\text{R}}(x)=
Ad + \sum_{j=1}^d (x_j - b_j)^2 - A \cos(2\pi(x_j - b_j)),\quad \text{with parameters } A\in\R,b\in\R^d
\end{align*}
and $b$  the location of the minimiser with $g^{\text{R}}(b)=0$. In order to highlight the efficiency of the componentwise variant of CBO we also slightly modify the Rastrigin function as follows:
\begin{align*}
g^{\hat{\text{R}}}(x)=\sum_{j=1}^d (x_j^2 - b_j^2)^2 - h \cdot \cos(k  x_j) + c \cdot x_j,\quad\text{with parameters } b\in\R^d, h,k,c\in\R.
\end{align*}
The sum-of-different-powers function is a separable function from the \texttt{bbob} test suite \cite{bbob2019} (therein referred to as f14), with the following definition $g^{\text{DP}}(x) = 
\sum_{j=1}^{d} |x_j|^{2 + 4\tfrac{j-1}{d-1}}$. 
Finally, in order to illustrate cases where the estimation stage does not succeed we also consider the squared $\ell^2$-norm $g^{\ell_2}(x) = \|x\|_2^2 = \sum_{j=1}^d x_j^2$.
It is clear that such a function cannot fulfil the identifiability condition from Assumption \ref{assump:identifiabilityHessianLinSepFunc}.

In Figure~\ref{fig:rotrast}, we first illustrate the behaviour for a simple two-dimensional setup with the Rastrigin function. 
Already here we can observe a failure mode of anisotropic CBO on $f = g(W \cdot )$ and how our proposed reparametrisation resolves it. 
Beyond that, we perform a quantitative benchmark which is displayed in Table~\ref{tab:success-rates}. 
Here, a run is counted as successful if
%
%\begin{align*}
$f(\hat{x}) - f(x^*) \leq 0.005$
%\end{align*}
%
is fulfilled,
where $\hat{x}$ denotes the best particle found in the algorithm and $x^*$ the best known minimum of $f$. 
The results in Table~\ref{tab:success-rates} demonstrate clearly the efficacy of our approach.
Importantly, despite the CBO schemes being allowed to use only $B - B_{H}$ queries due to the preceding estimation step, our approach still yields a drastic increase in performance. 
In particular, the splitting across dimensions can outperform the other CBO reparametrisation variants on $g^{\hat{\text{R}}}$. 
We highlight that here the number of queries is also shared across different dimensions, meaning each single one-dimensional CBO run can only use $\lfloor \tfrac{B - B_{H}}{d} \rfloor$ queries. Finally, the results in Table~\ref{tab:success-rates} show that CMA-ES can adapt to the underlying structure especially well for bowl-like functions, but struggles to do so for more irregular functions like Rastrigin. 

\begin{table}[t]
\centering%
\begin{tabular}{lcccc}
\toprule
Method & Rastrigin & Rastrigin-Like & Squared Norm & Different Powers \\
\midrule
CBO on $g$ with $\mathsf{D}^{\text{aniso}}$ & 100\% & 88\% & 100\% & 100\% \\
CBO on $f=g(W\cdot)$ with $\mathsf{D}^{\text{aniso}}$ & 0\% & 0\% & 0\% & 0\% \\
\midrule%
CBO on $\tilde{g} = f(\tilde{W}^\dagger \cdot)$ with $\mathsf{D}^{\text{aniso}}$ & 99\% (100\%) & 74\% (100\%) & 0\% (0\%) & 98\% (100\%) \\
CBO with $\mathsf{D}^{\text{transform}}$ & 99\% (100\%) & 75\% (100\%) & 0\% (0\%) & 99\% (100\%) \\
Componentwise CBO on $\tilde{g}$ & 98\% (100\%) & 99\% (100\%) & 1\% (0\%) & 100\% (100\%) \\
\midrule%
BIPOP-CMA-ES on $f=g(W \cdot)$ & 0\% & 0\% & 100\% & 100\% \\
\bottomrule
\end{tabular}
\caption{Success rates of the studied algorithms for different benchmark functions in $d=100$ averaged over $100$ runs. 
The percentage in brackets refers to the successfully recovered weights $w_j$, which is fulfilled if $\|\tilde{w}_j - w_j\|_2\leq 10^{-4}$. 
%The estimation relies on sufficient information provided by the Hessians of $f_j$. In the case of a quadratic functional $\phi(s)=s^2$, we observe that $\phi''(s)=2$ is a constant function and thus in this case we cannot expect to recover $W$, since the identifiability condition \eqref{eq:identifiabilityConditionHessianLinSep} is violated. 
}
\label{tab:success-rates}
\end{table}%
\section{Conclusion} \label{sec:conclusion}
In this work, we give a partial answer to an open question posed in \cite{bonandin2026exploiting} on how to make CBO exploit the structure of linearly separable functions. 
We leverage an efficient estimation method with provable recovery guarantees to find a suitable coordinate transformation, and we propose three different modifications to the CBO scheme employing this estimate. 
Numerical experiments clearly demonstrate that through this modification CBO is in fact able to exploit linearly separable functions efficiently, where standard CBO fails.
Since both stages of our method come with strong theoretical guarantees, this significantly expands the scope of efficient global nonconvex optimisation with convergence guarantees.

Our work opens up a multitude of interesting and promising questions for future research, and we highlight two: 
First, can we still obtain convergence guarantees as in \cite{bonandin2026exploiting} when the function is only approximately separable? 
Second, instead of a two-stage method, can we perform the estimation step online during the CBO run?

% \section{Thm commands}
% Here is an example:
% \begin{theo}\label{th:true}
%   Most theorems are true.
% \end{theo}

% \begin{proof}
%   Th. \ref{th:true} is obviously true.
% \end{proof}

% \begin{exam}\label{ex:good}
%   This should look like a good example.
% \end{exam}

% \begin{rema}
%   Can an example like Ex.~\ref{ex:good} give some insight in
%   Th.~\ref{th:true}'s proof?
% \end{rema}

% The next command determines the bibliography style. Please do not
% change this.
\bibliographystyle{plain+eid}

%  This inserts the bib file
\bibliography{refs}
%\bibliography{refs.bib}

\end{document}